\documentclass[a4paper,11pt]{article}

\usepackage[french]{babel}
\usepackage[utf8]{inputenc}
\usepackage[T1]{fontenc}
\usepackage{amsmath}
\usepackage{amsfonts}
\usepackage{amssymb}
\usepackage{dsfont}
\usepackage{graphicx}
\usepackage{hyperref}
\usepackage{xcolor}
\hypersetup{linkcolor=blue}
\hypersetup{colorlinks=true}
\usepackage[width=15cm, height=24cm]{geometry}

\usepackage{amsthm}
\newtheorem{thm}{Th\'eor\`eme}
\newtheorem{defi}{D\'efinition}
\newtheorem{lemme}{Lemme}

\title{Obstructions quadratiques à la contrôlabilité, \\ de la dimension finie à la dimension infinie}
\author{Frédéric Marbach}
\date{}

\newcommand{\ds}{\mathrm{d}s}
\newcommand{\dx}{\mathrm{d}x}
\newcommand{\dt}{\mathrm{d}t}

\newcommand{\Sone}{\mathcal{S}_1(0)}
\newcommand{\Stwo}{\mathcal{S}_2(0)}
\newcommand{\tm}{\mathcal{T}}
\newcommand{\M}{\mathcal{M}}
\newcommand{\norm}[2]{\left\| #1 \right\|_{#2}}
\newcommand{\ad}{\mathrm{ad}}
\newcommand{\R}{\mathbb{R}}
\newcommand{\N}{\mathbb{N}}
\newcommand{\C}{\mathbb{C}}

\begin{document}

\maketitle

\begin{abstract}
 On s'intéresse à la contrôlabilité locale en temps petit d'un système au voisinage
 d'un équilibre. \'Etant donné un petit temps imparti, une donnée initiale proche
 de l'équilibre et une donnée finale proche de l'équilibre, est-il possible de trouver
 un contrôle (un terme source) qui guide la solution du système depuis l'état initial
 vers l'état final souhaité dans le temps imparti ? La démarche naturelle consiste
 à commencer par étudier la contrôlabilité du système linéarisé au voisinage de 
 l'équilibre. Lorsque celui-ci n'est pas contrôlable, il est nécessaire de poursuivre
 le développement à l'ordre quadratique. 
 
 Dans cette note, on fait le lien entre plusieurs résultats récents autour de cette
 thématique, dans le cas particulier où le contrôle est scalaire.
 Ces résultants laissent penser que, dans ce cas, l'ordre quadratique ne peut qu'apporter des
 obstructions à la contrôlabilité. On commente notamment un résultat exhaustif obtenu en 
 collaboration avec Karine Beauchard en
 dimension finie et un résultat contenu dans la thèse de l'auteur qui fait
 apparaître un phénomène nouveau en dimension infinie.
\end{abstract}

\section{Introduction}

On présente ici la problématique dans le contexte des systèmes différentiels ordinaires
 pour simplifier l'exposition. 
La généralisation aux systèmes dont l'évolution est régie par des équations aux dérivées partielles est exposée dans les Sections~\ref{Section:pde.same} et~\ref{Section:pde.new}.

\subsection{Contrôlabilité locale en temps petit}

On note $x(t) \in \R^n$ l'état du système à l'instant $t \in [0,T]$. On suppose que l'évolution de cet état dépend du choix d'un contrôle $u(t) \in \R$ par l'opérateur. 
Ce choix correspond à l'idée d'un unique contrôle scalaire et constitue un cas 
particulier d'une situation plus générale où le contrôle prendrait ses valeurs
dans un espace de dimension supérieure. On suppose en outre que la dynamique est
déterminée par une fonction $f \in C^1(\R^n \times \R, \R^n)$ de sorte que :
\begin{equation} \label{fxu}
 \dot{x} = f(x,u).
\end{equation}
On s'intéresse à la contrôlabilité de~\eqref{fxu} au voisinage d'un équilibre.
Par translation, on peut supposer que cet équilibre est situé en $x = 0$ et $u = 0$
et que l'on a donc $f(0,0) = 0$. De multiples notions de contrôlabilité locale
existent. Dans cet exposé, on étudie la possiblité de ramener l'état à l'équilibre.
Plus précisément, on considère la notion suivante.

\begin{defi}
 On dit que le système~\eqref{fxu} est localement contrôlable en temps petit
 (vers l'équilibre nul) lorsque, pour tout temps $T > 0$, 
 pour tout taille maximale du contrôle $\eta > 0$, il existe un voisinage de
 l'origine de taille $\delta > 0$ tel que, pour toute donnée initiale $x^*$
 avec $|x^*| \leqslant \delta$, il existe un contrôle $u \in L^\infty((0,T),\R)$
 de norme plus petite que $\eta$ telle que la solution de~\eqref{fxu}
 associée à ce contrôle et à la donnée initiale $x(0)=x^*$ vérifie $x(T) = 0$.
\end{defi}

\subsection{Linéarisation, critère de Kalman et test linéaire}

La démarche naturelle pour étudier la contrôlabilité locale d'un système consiste
à linéariser sa dynamique au voisinage de l'équilibre. Ainsi, on regarde le système :
\begin{equation} \label{def.Ab}
 \dot{y} = A y + u b
 \quad \textrm{où} \quad
 A := \frac{\partial f}{\partial x}(0,0) \in \mathcal{M}_n(\R)
 \quad \textrm{et} \quad
 b := \frac{\partial f}{\partial u}(0,0) \in \R^n.
\end{equation}
Il est bien connu (voir~\cite[Théorème 10]{MR0155070}) qu'un tel système linéaire 
est contrôlable si et seulement s'il vérifie le critère de Kalman :
\begin{equation} \label{kalman}
 \mathrm{rang~} \{ b, Ab, A^2b, \ldots A^{n-1}b \} = n.
\end{equation}

\begin{thm}[voir {\cite{MR889459}}] \label{thm1}
 Soit $f \in C^1(\R^n\times\R,\R^n)$. On suppose que la paire $(A,b)$
 définie en~\eqref{def.Ab} vérifie~\eqref{kalman}. 
 Alors le système~\eqref{fxu}
 est localement contrôlable en temps petit.
\end{thm}

Le test linéaire fournit donc une condition suffisante de contrôlabilité locale
en temps petit pour un système. 
Cette condition n'est évidemment pas nécessaire, comme on peut le vérifier à partir
de l'exemple suivant :
\begin{equation}
 \left\{
  \begin{aligned}
   \dot{x}_1 & = u, \\
   \dot{x}_2 & = x_1^3.
  \end{aligned}
 \right.
\end{equation}  
Pour de tels systèmes, il faut poursuivre le développement plus loin pour conclure.

\subsection{Systèmes affines en le contrôle et crochets de Lie}

Une sous-classe importante des systèmes non-linéaires généraux est celle des 
systèmes dont la dynamique dépend du contrôle de manière affine  :
\begin{equation} \label{f0f1}
 \dot{x} = f_0(x) + u f_1(x).
\end{equation}
Ces systèmes sont pertinents du point de vue des applications 
(car il est fréquent que les contrôles soient des forces et interviennent donc
de manière affine dans la dynamique) mais aussi du point de vue de l'étude 
mathématique (à la fois comme première étape, comme développement de Taylor en $u \ll 1$
à l'ordre un, ou même par réduction de n'importe quel système sous cette forme
après une transformation de type intégrateur).

Pourtant, même pour ces systèmes, on ne connait pas de condition nécessaire et
suffisante de contrôlabilité locale en temps petit. Si $f_0$ et $f_1$ 
sont des champs de vecteurs analytiques, la condition suivante est nécessaire 
(voir~\cite{MR0149402}) :
\begin{equation} \label{lie.rank}
 \mathrm{Lie~} \{f_0, f_1\}(0) = \R^n.
\end{equation}
Si~\eqref{lie.rank} ne tient pas, alors l'état est contraint d'évoluer dans
une sous-variété stricte.
De plus, pour un système lisse vérifiant~\eqref{lie.rank}, si l'on définit, 
pour $k\in \N$, $\mathcal{S}_k(0)$ comme étant l'espace vectoriel
engendré par tous les crochets de Lie itérés de $f_0$ et de $f_1$ contenant
$f_1$ au plus $k$ fois évalués en l'équilibre, la condition suivante est suffisante 
(voir~\cite{MR710995}) :
\begin{equation} \label{even.odd}
 \forall k \in \N, \quad \mathcal{S}_{2k+2} (0) = \mathcal{S}_{2k+1}(0). 
\end{equation}
L'interprétation habituelle de~\eqref{even.odd} est que les crochets
de type pair (en $f_1$) sont ``mauvais'' (car associés à un signe)
et doivent être compensés par de ``bons'' crochets de type impair.
L'étude des bons et des mauvais crochets a occupé de nombreux auteurs
(voir~\cite{MR1061394}). 

\subsection{Les premières obstructions quadratiques}

Notons que l'espace $\mathcal{S}_1(0)$ correspond à celui qui
intervient à l'ordre linéaire dans~\eqref{kalman}. Dans cette
note, on s'intéresse à la condition $\mathcal{S}_2(0) = \mathcal{S}_1(0)$
car l'espace $\mathcal{S}_2(0)$ est associé aux propriétés des 
approximations quadratiques.
La première violation possible de la condition $\mathcal{S}_2(0) = \mathcal{S}_1(0)$
se produit lorsque:
\begin{equation} \label{sussmann1}
 [f_1, [f_0, f_1]](0) \notin \mathcal{S}_1(0).
\end{equation}
Sous cette hypothèse, Sussmann démontre le résultat suivant.

\begin{thm}[voir \cite{MR710995}]
 Soit $f_0$ et $f_1$ deux champs de vecteurs analytiques. On suppose~\eqref{sussmann1}.
 Alors le système~\eqref{f0f1} n'est pas localement contrôlable en temps petit.
\end{thm}

La violation quadratique suivante se produit lorsque :
\begin{equation} \label{sussmann2}
 [[f_0,f_1], [f_0, [f_0, f_1]]](0) \notin \mathcal{S}_1(0).
\end{equation}
Or, la condition~\eqref{sussmann2} n'empêche pas un système d'être localement
contrôlable en temps petit, comme le montre le système suivant, étudié par
Sussmann dans~\cite{MR710995}:
\begin{equation}
 \left\{
 \begin{aligned}
 \dot{x}_1 & = u, \\
 \dot{x}_2 & = x_1, \\
 \dot{x}_3 & = x_2^2 + x_1^3.
 \end{aligned}
 \right.
\end{equation}
Ainsi, l'étude des approximations quadratiques avait un peu été mise de côté,
car cet exemple semblait aller à l'encontre de l'intuition que la condition
$\mathcal{S}_2(0) = \mathcal{S}_1(0)$ devait être nécessaire en un certain sens,
éventuellement plus faible.

\section{Classification des obstructions en dimension finie}

Dans \cite{2017arXiv170507991B}, on classifie les obstructions quadratiques à la contrôlabilité
locale en temps petit pour des systèmes mono-contrôle en dimension finie. 
Le résultat principal prend la forme d'une alternative quadratique.
Lorsque $\mathcal{S}_2(0) \not \subset \mathcal{S}_1(0)$, on démontre
l'existence d'une dérive quadratique signée. Lorsque
$\mathcal{S}_2(0) = \mathcal{S}_1(0)$, on démontre que les termes
quadratiques n'augmentent pas la dimension de l'espace atteignable
par rapport à l'ordre linéaire.

\subsection{Notations et définitions}

Une des idées clefs de ce travail est d'introduire la notion suivante, qui met en valeur 
le rôle joué par l'hypothèse de norme sur les contrôles.

\begin{defi}
 Soit $(E_T, \norm{\cdot}{E_T})_{T>0}$ une famille d'espaces vectoriels normés
 de fonctions scalaires définies sur $[0,T]$.
 On dit que le système~\eqref{fxu} est $E$ localement contrôlable en temps petit
 (vers l'équilibre nul) lorsque, pour tout temps $T > 0$, 
 pour tout taille maximale du contrôle $\eta > 0$, il existe un voisinage de
 l'origine de taille $\delta > 0$ tel que, pour toute donnée initiale $x^*$
 avec $|x^*| \leqslant \delta$, il existe un contrôle $u \in L^\infty((0,T),\R) \cap E_T$
 de norme dans $E_T$ plus petite que $\eta$ telle que la solution de~\eqref{fxu}
 associée à ce contrôle et à la donnée initiale $x(0)=x^*$ vérifie $x(T) = 0$.
\end{defi}

Par ailleurs, pour deux champs de vecteurs $f_0$ et $f_1$, on définit par récurrence pour $k \in \N$ les champs
de vecteurs $\ad^k_{f_0}(f_1)$ comme étant le crochet de Lie itéré à gauche où $f_0$ 
intervient $k$ fois et $f_1$ une seule fois.

Enfin, considérant le produit scalaire usuel sur $\R^n$, on 
note $\mathbb{P}$ la projection orthogonale sur le sous-espace vectoriel $\mathcal{S}_1(0)$ et $\mathbb{P}^\perp := \mathrm{Id}
- \mathbb{P}$.

\subsection{\'Enoncé du résultat}

On énonce ci-dessous le résultat dans sa forme la plus simple pour les
systèmes affines en le contrôle. Nous renvoyons le lecteur à \cite{2017arXiv170507991B} pour
une version du théorème incluant le cas des systèmes non-linéaires
quelconques.

\begin{thm}
 Soit $f_0, f_1 \in C^{\infty}(\R^n,\R^n)$ avec $f_0(0) = 0$
 et $d:=\mathrm{dim~}\Sone$. Il existe une fonction $G \in C^\infty(\Sone,\Sone^\perp)$ vérifiant $G(0) = 0$ et $G'(0) = 0$, 
 telle que l'alternative suivante tienne :
 \begin{itemize}
  
  \item \textbf{Lorsque} $\Stwo = \Sone$, à un terme d'erreur 
  cubique en le contrôle près, l'état vit dans une variété lisse
  $\M \subset \R^n$ de dimension $d$ donnée par le graphe de $G$ :
  \begin{equation} \label{def.MG1}
  \M := \left\{ p + G(p) ; \enskip p \in \Sone \right\}.
  \end{equation} 
  Plus précisément, pour tout $\tm > 0$, il existe $C, \eta > 0$ tels que, 
  pour tout $T \in (0,\tm]$, pour toute trajectoire 
  $(x,u) \in C^0([0,T],\R^n) \times L^1((0,T),\R)$ de \eqref{f0f1} 
  avec $x(0) = 0$ et $\norm{u}{W^{-1,\infty}} \leqslant \eta$, on ait :
  \begin{equation}
  \forall t \in [0,T], \quad
  \left| \mathbb{P}^\perp x(t) - G \left( \mathbb{P} x(t) \right) \right| 		
  \leqslant C \norm{u}{W^{-1,3}}^3.
  \end{equation}
  
  \item \textbf{Lorsque} $\Stwo \not \subset \Sone$, pour des temps
  assez petits et des contrôle assez réguliers, l'état dérive 
  par rapport à la variété invariante $\M$ dans une direction donnée.
  Plus précisément :
  \begin{itemize}
   
   \item Il existe $k \in \{1, \ldots, d \}$ tel que :
   \begin{align} 
   \label{hyp_j}
   [ \ad_{f_0}^{j-1}(f_1) , \ad_{f_0}^{j}(f_1) ] (0) & \in \Sone 
   \quad \text{ pour } 1 \leqslant j < k, \\
   \label{hyp_k}
   [ \ad_{f_0}^{k-1}(f_1) , \ad_{f_0}^{k}(f_1) ] (0) & \notin \Sone.
   \end{align}
   
   \item Le système~\eqref{f0f1} n'est pas $W^{2k-3,\infty}$
   localement contrôlable en temps petit.
   
   \item Il existe $T^*>0$ tel que, 
   pour tout $\tm \in (0,T^*)$, il existe $\eta > 0$ tel que,
   pour tout $T \in (0,\tm]$ et toute trajectoire 
   $(x,u) \in C^0([0,T],\R^n) \times L^1((0,T),\R)$
   de~\eqref{f0f1} avec $x(0) = 0$ 
   et $u \in W^{2k-3,\infty}_0$ 
   et $\norm{u}{W^{2k-3,\infty}} \leqslant \eta$, on ait :
   \begin{equation} \label{drift}
   \forall t \in [0,T], \quad \left\langle \mathbb{P}^\perp x(t)- 
   G \left( \mathbb{P} x(t) \right) , d_k \right\rangle \geqslant 0,
   \end{equation}
   où la direction de dérive $d_k \neq 0$ est donnée par :
   \begin{equation} \label{def.dk}
   d_k := - \mathbb{P}^\perp 
   [ \ad_{f_0}^{k-1}(f_1) , \ad_{f_0}^{k}(f_1) ] (0).
   \end{equation}
  \end{itemize}
 \end{itemize}
\end{thm}

\subsection{Commentaires}

\begin{itemize}
 \item Le temps critique $T^*$ est caractérisé par une propriété de coercivité d'une forme quadratique correspondant à une approximation d'ordre deux du système.
 \item Pour $\tm$ assez petit (ou pour tout $\tm < T^*$ pour certains systèmes), l'inégalité~\eqref{drift} tient en réalité dans un sens plus fort :
 \begin{equation} \label{drift.strong}
 \forall t \in [0,T], \quad \left\langle \mathbb{P}^\perp x(t)- 
 G \left( \mathbb{P} x(t) \right) , d_k \right\rangle \geqslant C(\tm) \int_0^t u_k^2(s) \ds,
 \end{equation}
 où $u_k$ est la $k$-ième primitive itérée du contrôle $u$.
 \item Il est possible de construire des systèmes pour lesquels les espaces fonctionnels dans lesquels sont exprimés les hypothèses de petitesse pour les contrôles sont optimaux (au moins au sein des espaces de Sobolev) et ce même au sein de la classe des systèmes bilinéaires.
 De même, la coercivité obtenue dans~\eqref{drift.strong} est optimale.
 \item L'inégalité~\eqref{drift} tient pour tout $t \in [0,T]$ uniquement si 
 $u \in W^{2k-3}_0$, ce qui signifie pour nous que $u(0) = u'(0) = ... = u^{(2k-2)}(0) = 0$. Pour un contrôle quelconque, l'inégalité ne tient pas dès $t = 0$ mais tient tout de 
 même au temps final, ce qui permet de nier la contrôlabilité locale en temps
 petit au sens $W^{2k-3}$.
 \item Une façon d'interpréter ce résultat est de dire que la condition $\mathcal{S}_2(0) = \mathcal{S}_1(0)$ est nécessaire pour qu'un système soit localement contrôlable en temps petit avec des contrôles réguliers. En effet, si $\mathcal{S}_2(0) \neq  \mathcal{S}_1(0)$, le théorème indique qu'il existe une régularité maximale critique pour les contrôles : ils ne peuvent pas être petits dans des espaces plus réguliers.
 \item Le théorème est énoncé par simplicité avec des champs $f_0$ et $f_1$ lisses.
 Cependant, toutes les conclusions sont préservées si on prend des champs seulement $C^3$. En effet, on peut vérifier que l'ensemble des objets en jeu sont encore
 bien définis et que toute la preuve tient.
\end{itemize}

\subsection{\'Eléments de démonstration}

La démonstration repose sur plusieurs idées clefs :
\begin{itemize}
 \item Tout d'abord, on se ramène à l'étude de systèmes vérifiant $A^d = 0$,
 où la matrice $A$ est définie par~\eqref{def.Ab} et où $d$ est la dimension
 de l'espace contrôlable. Cette simplification se fait par une transformation de type Brunovsk\'y
 et permet de réduire la dynamique linéaire à un intégrateur simple.
 \item Ensuite, on introduit $d$ systèmes auxiliaires en composant successivement la dynamique par les flots $\phi_j$ pour $1\leqslant j\leqslant d$ correspondants aux champs de vecteurs $(-1)^{j-1}\ad^{j-1}_{f_0}(f_1)$. Ces transformations permettent de prendre en compte entièrement la dynamique linéaire et de relever tous les termes qui vont constituer la variété invariante de référence.
 \item On étudie alors l'approximation quadratique du $d$-ième système auxiliaire et l'on montre qu'elle fait intervenir des directions de dérive données par~\eqref{def.dk} et pondérées par des quantités qui ont un signe en temps petit et sont même 
 coercives pour des normes de Sobolev $H^{-k}$ du contrôle.
 \item Enfin, on conclut que la dérive observée à l'ordre quadratique persiste sur le système initial en utilisant des inégalités d'interpolation du type Gagliardo-Nirenberg (voir~\cite{MR0109940}) afin d'absorber les termes cubiques (ou d'ordre supérieur) par le produit de la dérive coercive (dans une norme faible) et une norme forte sur le contrôle (dont on fait l'hypothèse qu'elle est suffisamment petite).
\end{itemize}

\section{Situations analogues en dimension infinie} \label{Section:pde.same}

En dimension infinie, pour étudier les propriétés de contrôlabilité locale d'un
système, la démarche présentée en introduction reste valable : on commence toujours par considérer le système linéarisé. Lorsque celui-ci est contrôlable, le passage au non-linéaire se fait aussi par une approche de type point-fixe (même s'il n'existe pas de résultat aussi général que le Théorème~\ref{thm1}). En effet, l'étude du linéarisé
et le passage au non-linéaire sont plus difficiles et doivent être adaptés pour chaque type d'EDP considéré.

Lorsque le linéarisé ne permet pas de conclure, plusieurs études de cas particuliers ont déjà été menées. On constate qu'on peut retrouver les comportements trouvés en
dimension finie : le long des directions perdues à l'ordre linéaire, l'ordre quadratique apporte parfois une dérive signée et parfois un comportement de type variété invariante. Donnons, à titre d'illustration, deux exemples opposés.

\subsection{Une équation de Schrödinger}

On considère le système contrôlé suivant :
\begin{equation} \label{schro}
 \left\{
  \begin{aligned}
   i \psi_t + \psi_{xx} + u(t) \mu(x) \psi & = 0, \\
   \psi(t,0) = \psi(t,1) & = 0.
  \end{aligned}
 \right.
\end{equation}
Ce système représente une particule quantique, coincée dans un puits de potentiel infini et soumise à un champ électrique uniforme d'amplitude $u(t)$. Le paramètre $\mu \in H^3((0,1),\R)$ décrit le moment dipolaire de la particule. La norme $L^2((0,1),\C)$ de la fonction d'onde $\psi$ est invariante pendant l'évolution. On étudie la contrôlabilité
de ce système au voisinage de l'état fondamental $\psi_1(t,x) := \varphi_1(x) e^{-i\pi^2 t}$, où $\varphi_1(x) := \sqrt{2} \sin (\pi x)$.

Il est connu (voir~\cite{MR2732927}) que ce système est localement contrôlable en temps petit lorsque le moment $\mu$ est suffisamment riche au sens où il existe $c > 0$ tel que :
\begin{equation}
 \forall k \in \N^*, \quad \left|\langle\mu\varphi_1,\varphi_k\rangle\right| \geqslant\frac{c}{k^3}.
\end{equation}
A contrario, si une des projections est nulle, la direction correspondante est
perdue à l'ordre linéaire et il faut étudier l'approximation quadratique du système
pour pouvoir conclure.
\begin{thm}[voir \cite{MR3167929}]
 Soit $\mu \in H^3((0,1),\R)$. On suppose qu'il existe $k \in \N^*$ tel que :
 \begin{equation} \label{eq.ak}
  \langle \mu \varphi_1, \varphi_k \rangle = 0 \quad \textrm{et} \quad 
  \alpha_k := \langle (\mu')^2 \varphi_1, \varphi_k \rangle \neq 0.
 \end{equation}
 Alors il existe $T > 0$, $\eta > 0$ et une suite d'états finaux $\psi^\delta$
 arbitrairement proches de $\psi_1(T)$
 tels que, pour tout contrôle
 $u \in L^2((0,T),\R)$ avec $\norm{u}{L^2} \leqslant \eta$, la solution de~\eqref{schro}
 partant de l'état initial $\psi_1(0,x)$ ne vérifie pas $\psi(T) = \psi^\delta$.
\end{thm}
Ainsi, le système n'est pas localement contrôlable en temps petit.
La démonstration repose sur l'étude d'un système auxiliaire et l'établissement d'une dérive quadratique quantifiée par la norme $H^{-1}(0,T)$ du contrôle. Les auteurs font
 remarquer que la condition $\alpha_k \neq 0$ est indispensable à leur argument.
 
\`A la lumière des travaux en dimension finie, on peut maintenant expliquer plus précisément cette condition. En effet, supposer $\alpha_k \neq 0$ revient à supposer
que l'analogue du premier crochet de Lie $[f_1, [f_0, f_1]](0)$ n'appartient pas à $\mathcal{S}_1(0)$.
Dans ce contexte, $f_0$ est l'opérateur $i \partial_{xx}$ et $f_1$ la multiplication par $i \mu$. Ainsi $[f_1, [f_0, f_1]]$ est la multiplication par $2i(\mu')^2$.
Comme on effectue une étude locale au voisinage de l'état de base $\varphi_1$, 
la direction du crochet est $2i(\mu')^2\varphi_1$ et l'hypothèse~\eqref{eq.ak}
traduit exactement l'idée que ce crochet a une composante non nulle en dehors
de l'espace contrôlable. Il est donc normal d'observer une dérive du premier type.
De plus, l'hypothèse de petitesse du contrôle obtenue dans l'article
peut être ramenée à $H^{-1}$, comme c'est le cas en dimension finie.

Enfin, dans le cas où $\alpha_k = 0$, il faut vraisemblablement poursuivre 
le développement plus loin pour étudier la présence ou l'absence d'obstruction
quadratiques liés à des crochets d'ordre supérieur. Bien sûr, cette étude
se complique fortement en dimension infinie pour des questions de régularité :
il faut que $\mu$ soit suffisamment régulier pour pouvoir définir
les crochets d'ordre supérieur.
 
\subsection{Une équation de Korteweg-de-Vries}

Soit $L > 0$. On considère le système contrôlé suivant :
\begin{equation} \label{kdv}
 \left\{
  \begin{aligned}
   \psi_t + \psi_x + \psi_{xxx} + \psi \psi_x & = 0, \\
   \psi(t,0) = \psi(t,L) & = 0, \\
   \psi_x(T,L) & = 0.
  \end{aligned}
 \right.
\end{equation}
Ce système modélise la propagation de vagues de petite amplitude à la surface 
d'un canal uniforme. Il est connu que ce système est localement contrôlable au voisinage de l'équilibre
nul lorsque la longueur $L$ ne fait pas partie d'un ensemble discret de longueurs dites
critiques (voir~\cite{MR1440078}). Pour ces longueurs non critiques, la contrôlabilité s'obtient à partir de l'étude du système linéarisé, puis par point fixe. Pour les longueurs critiques en revanche, l'étude du système linéarisé ne permet pas de conclure et la contrôlabilité locale en temps petit n'est pas encore connue dans tous les cas.
Pour la première longueur critique $L = 2\pi$, on a le résultat suivant :

\begin{thm}[voir \cite{MR2060480}]
 Soit $L = 2\pi$, $T > 0$ et $\eta > 0$. Il existe $\delta > 0$ tel que, pour
 toutes données initiales et finales $\psi^*, \psi^\dagger \in L^2(0,L)$
 de taille plus petite que $\delta$, il existe $u \in L^2(0,T)$ de norme
 plus petite que $\eta$ telle que la trajectoire 
 de~\eqref{kdv} partant de $\psi^*$ arrive en $\psi^\dagger$.
\end{thm}
 
Dans la formulation initiale de ce résultat, le contrôle n'apparaît pas mais est vu comme la trace appropriée de la trajectoire : $\psi_x(t,L)$. En considérant l'espace fonctionnel dans lequel vivent les trajectoires construites, on en déduit cependant la formulation ci-dessus.

Pour démontrer ce résultat, les auteurs effectuent un développement à l'ordre 3 de la dynamique. En effet, à l'ordre linéaire, la direction $(1-\cos(x))$ est perdue.
C'est la seule direction perdue et il faut poursuivre le développement pour comprendre si la présence du terme non-linéaire $\psi \psi_x$ permet ou non de la récupérer. Les auteurs démontrent qu'il est possible de bouger dans les deux sens $\pm (1-\cos(x))$
à l'aide des termes cubiques.

Or, ils démontrent aussi que, bien que cela aurait été plus simple, il n'est pas possible d'obtenir leur résultat en s'arrêtant à l'ordre deux. En effet, ils démontrent que, si l'on considère une trajectoire $y$ de :
\begin{equation} \label{kdv1}
\left\{
\begin{aligned}
y_t + y_x + y_{xxx} & = 0, \\
y(t,0) = y(t,L) & = 0, \\
y(0,x) = y(T,x) & = 0
\end{aligned}
\right.
\end{equation}
et une trajectoire $z$ de :
\begin{equation} \label{kdv2}
\left\{
\begin{aligned}
z_t + z_x + z_{xxx} & = - y y_x, \\
z(t,0) = z(t,L) & = 0, \\
z(0,x) & = 0,
\end{aligned}
\right.
\end{equation}
alors l'ordre quadratique vérifie :
\begin{equation}
 \int_0^L z(T,x) (1-\cos(x)) \dx = 0.
\end{equation}
Ceci fait penser au cas ``variété invariante'' du théorème en dimension finie :
prescrire le mouvement de l'ordre linéaire impose la position de l'ordre
quadratique dans le cas où $\mathcal{S}_2(0) = \mathcal{S}_1(0)$. 
Bien sûr, il n'est pas facile de donner un sens à cette égalité dans 
ce contexte assez singulier où le contrôle agit au bord.

\section{Un phénomène nouveau en dimension infinie} \label{Section:pde.new}

\subsection{\'Enoncé du résultat}

En dimension infinie, de nouvelles obstructions quadratiques peuvent apparaître,
comme en témoigne le résultat suivant pour une équation de Burgers, 
dont la démonstration fait intervenir une dérive d'un nouveau type, caractérisée
par la norme du contrôle dans~$H^{-5/4}$.
Pour $T > 0$ et $\psi_0 \in H^1_0(0,1)$ une donnée initiale, on considère le système 
de contrôle suivant, impliquant un unique contrôle scalaire $u : [0,T] \to \R$ :
\begin{equation} \label{OQB.2}
\left\{
\begin{aligned}
\psi_t + \psi \psi_x - \psi_{xx} & = u(t)
& \quad \quad [0,T] \times [0,1], \\
\psi(t, 0) & = 0 & [0,T], \\
\psi(t, 1) & = 0 & [0,T], \\		
\psi(0, x) & = \psi_0(x) & [0,1].
\end{aligned}
\right.
\end{equation}

\begin{thm}[voir \cite{2015arXiv151104995M}] \label{intro_thm.loc.u}
 Le système~\eqref{OQB.2} n'est pas localement contrôlable à zéro en temps petit.
 Pour tout $\eta > 0$, il existe $T > 0$ tel que, pour tout $\delta > 0$, il 
 existe une donnée initiale $\psi_0 \in H^1_0(0,1)$ avec $|\psi_0|_{H^1_0} \leqslant 
 \delta$	telle que, pour tout contrôle $u \in L^2(0,T)$ satisfaisant 
 $|u|_{L^2} \leqslant \eta$, la solution $\psi$ de~\eqref{OQB.2} 
 est telle que $\psi(T,\cdot)\neq 0$. De plus, la donnée 
 initiale $\psi_0$ peut être choisie très régulière (et même polynomiale en $x$). 
\end{thm}

\subsection{\'Eléments de démonstration}

On considère un temps imparti $T \ll 1$ que l'on note $\varepsilon$ pour 
illustrer la limite envisagée. On effectue le changement 
d'échelle habituel $\psi(t,x) \leftarrow \varepsilon \psi(\varepsilon t, x)$
et $u(t) \leftarrow \varepsilon^2 u(\varepsilon t)$
pour se ramener à une étude pour $t 
\in [0,1]$ du système :
\begin{equation} \label{OQB.3}
\left\{
\begin{aligned}
{\psi}_t + {\psi} {\psi}_x 
- \varepsilon {\psi}_{xx} & = u(t)
& \quad \quad [0,1] \times [0,1], \\
{\psi}(t, 0) & = 0 & [0,1], \\
{\psi}(t, 1) & = 0 & [0,1], \\		
{\psi}(0, x) & = \varepsilon \psi_0(x) & [0,1].
\end{aligned}
\right.
\end{equation}
On commence par étudier le cas $\psi_0 = 0$. Le cas général $\psi_0 \neq 0$ s'en déduit
en décomposant la solution complète comme étant la somme de la solution partant de $\varepsilon \psi_0$
avec un contrôle nul, de la solution avec donnée initiale nulle mais un contrôle $u$
et d'un terme de reste petit devant chacune des deux. 
La première étape est de 
considérer le système linéarisé de~\eqref{OQB.3} autour de $\psi = 0$. 
En notant $y$ l'état linéarisé, on obtient :
\begin{equation} \label{OQB.4}
\left\{
\begin{aligned}
y_t - \varepsilon y_{xx} & = u(t)
& \quad \quad [0,1] \times [0,1], \\
y(t, 0) & = 0 & [0,1], \\
y(t, 1) & = 0 & [0,1], \\
y(0, x) & = 0 & [0,1].
\end{aligned}
\right.
\end{equation}
Or, le système~\eqref{OQB.4} n'est pas contrôlable. 
En effet, le second membre $u(t)$ est à voir comme $u(t) \chi_{[0,1]}(x)$ où 
la fonction caractéristique $\chi_{[0,1]}$ est une fonction paire sur ce même 
segment (par rapport à $1/2$). Par conséquent, le contrôle n'agit que sur les modes pairs du système. 
Les modes impairs évoluent librement. On peut démontrer 
que le système~\eqref{OQB.4}, restreint aux modes pairs, est 
contrôlable. Ceci se fait, par exemple, en résolvant un problème de moments 
comme dans~\cite{MR0335014}.

On développe à l'ordre suivant en écrivant formellement 
$\psi = \eta y + \eta^2 z + \mathcal{O}(\eta^3)$. On obtient l'équation pour le 
second ordre : 
\begin{equation} \label{OQB.5}
\left\{
\begin{aligned}
z_t - \varepsilon z_{xx} & = - y y_x
& \quad \quad [0,1] \times [0,1], \\
z(t, 0) & = 0 & [0,1], \\
z(t, 1) & = 0 & [0,1], \\		
z(0, x) & = 0 & [0,1].
\end{aligned}
\right.
\end{equation}
L'intuition est que la quantité quadratique présente au second membre 
de~\eqref{OQB.5} va permettre, via un argument de signe, de nier la 
possibilité de retourner à zéro. On cherche une 
projection de $z$ dont on peut démontrer qu'elle 
a un signe en temps petit. Soit $\rho \in H^1_0(0,1)$ un profil de projection à 
fixer plus tard. Puisque $y$ dépend linéairement de $u$ et $z$ quadratiquement 
de $y$, il est raisonnable de penser qu'on peut écrire au temps final $t = 1$ 
une formule du type :
\begin{equation} \label{OQB.6}
\int_0^1 z(1, x) \rho(x) \dx 
= \int_0^1 \int_0^1 K^\varepsilon(s_1,s_2) u(s_1) u(s_2) \ds_1 \ds_2.
\end{equation}
Un argument de signe sur la projection de $z$ se traduit par une éventuelle 
coercivité de l'opérateur intégral correspondant au noyau $K^\varepsilon$. En 
considérant les systèmes adjoints de~\eqref{OQB.4} 
et~\eqref{OQB.5}, il est possible de démontrer 
l'égalité~\eqref{OQB.6} et de donner l'expression du noyau :
\begin{equation} \label{OQB.7}
K^\varepsilon(s_1, s_2) 
= \frac{1}{2} \int_{s_1 \vee s_2}^1 \int_0^1 \Phi^\varepsilon_x(1-t,x) 
G^\varepsilon(t-s_1, x) G^\varepsilon(t-s_2, x) \dx \dt,
\end{equation}
où $\Phi^\varepsilon$ et $G^\varepsilon$ sont données comme les 
solutions des systèmes de type ``chaleur'' suivants :
\begin{equation} \label{intro_eq.thm.loc.u.G}
\left\{ 
\begin{aligned}
G^\varepsilon_t - \varepsilon G^\varepsilon_{xx} & = 0 & \quad 
[0,1] \times [0,1], \\
G^\varepsilon(t,0) &= 0 & [0,1], \\
G^\varepsilon(t,1) &= 0 & [0,1], \\
G^\varepsilon(0,x) &= 1 & [0,1]~
\end{aligned}
\right.
\end{equation}
et
\begin{equation} \label{intro_eq.thm.loc.u.Phi}
\left\{ 
\begin{aligned}
\Phi^\varepsilon_t - \varepsilon \Phi^\varepsilon_{xx} & = 0 & 
\quad [0,1] \times [0,1], \\
\Phi^\varepsilon(t,0) &= 0 & [0,1], \\
\Phi^\varepsilon(t,1) &= 0 & [0,1], \\
\Phi^\varepsilon(0,x) &= \rho(x) & [0,1].
\end{aligned}
\right.
\end{equation}
Comme on raisonne en temps petit, on commence par s'intéresser au noyau limite 
lorsque $\varepsilon \rightarrow 0$. Sous réserve de 
quelques hypothèses correspondant à un choix approprié de $\rho$, 
on obtient l'équivalent suivant, valable ponctuellement en $s_1, s_2$ :
\begin{equation} \label{intro_eq.thm.loc.u.5}
K^\varepsilon(s_1,s_2) \sim 
\sqrt{\varepsilon} \left(|2 - s_1 - s_2|^{3/2}-|s_1-s_2|^{3/2}\right).
\end{equation}
On est amené à étudier le noyau limite du second membre 
de~\eqref{intro_eq.thm.loc.u.5} :
\begin{equation} \label{intro_eq.thm.loc.u.6}
\begin{split}
\langle K^0 u, u \rangle 
&:=\int_0^1\int_0^1 
\left(|2-s_1-s_2|^{3/2}-|s_1-s_2|^{3/2}\right)u(s_1)u(s_2)\ds_1\ds_2\\
&=\frac{3}{4}\int_0^1\int_0^1 
\left(|s_1-s_2|^{-1/2}+|2-s_1-s_2|^{-1/2}\right)U(s_1)U(s_2)\ds_1\ds_2,
\end{split}
\end{equation}
où l'on a introduit $U$ la primitive de $u$ s'annulant en $0$. On vérifie 
que le premier terme définit la norme $H^{-1/4}$ de $U$ 
et que le second est positif. Ainsi, on conjecture que :
\begin{equation} \label{intro_eq.thm.loc.u.7}
\langle K^\varepsilon u, u \rangle \geqslant \sqrt{\varepsilon} 
|u|_{H^{-5/4}}^2.
\end{equation}
L'établissement rigoureux de~\eqref{intro_eq.thm.loc.u.7} demande en réalité beaucoup 
de travail car l'asymptotique~\eqref{intro_eq.thm.loc.u.5} n'est que formel. Il 
faut vérifier que les noyaux résiduels sont petits (en~$\varepsilon$) et 
réguliers (au sens des espaces de fonctions sur lesquels ils sont continus).

On exprime la différence 
$R^\varepsilon(s_1,s_2) := K^\varepsilon(s_1,s_2) - \sqrt{\varepsilon} 
K^0(s_1,s_2)$ comme une fonction de $s_1$ et $s_2$. Il faut alors disposer d'un 
critère quantitatif, vérifiable sur un noyau quelconque $L(s_1,s_2)$, permettant de 
démontrer l'existence d'une constante $C(L) \geqslant 0$ 
telle que $|\langle Lu, u \rangle| \leqslant C(L)
|u|_{H^{-5/4}}^2$ et la constante $C(L)$ doit être calculable (ou au moins 
estimable) à partir de quantités dépendant relativement explicitement de $L$. 
On démontre ainsi que $|\langle R^\varepsilon u, u \rangle| \leqslant C 
\varepsilon^{3/2} |u|_{H^{-5/4}}^2$. Pour cela, on commence car calculer
explicitement $R^\varepsilon(s_1,s_2)$. Puis, on effectue la même
manipulation d'intégration par parties que dans~\eqref{intro_eq.thm.loc.u.6}.
Ainsi, on se ramène a étudier $\partial_{12} R^\varepsilon$ et on cherche
à démontrer
une estimation de la forme $|\langle \partial_{12} R^\varepsilon U, U \rangle| 
\leqslant C \varepsilon^{3/2} |U|_{H^{-1/4}}^2$. Une étape clef est l'utilisation
de la théorie des opérateurs intégraux faiblement singuliers, notamment
développée dans~\cite{MR1048075, MR1397491}. En particulier, on dispose du
critère suivant :

\begin{lemme} \label{lemma.wsio}
 Soit $L$ une fonction continue définie sur $\left\{ (t,s) \in (0,1) 
 \times (0,1) \textrm{, t.q. } t \neq s \right\}$. On suppose qu'il existe 
 $\kappa > 0$ et $\frac{1}{2} < \delta \leq 1$, tels que :
 \begin{align}
 \label{est.L.1}
 |L(t,s)| & \leq \kappa |t-s|^{-\frac{1}{2}}, \\
 \label{est.L.2}
 |L(t,s)-L(t',s)| & \leq \kappa |t-t'|^\delta |t-s|^{-\frac{1}{2}-\delta},
 \quad \textrm{pour } |t-t'|\leq\frac{1}{2}|t-s|, \\
 \label{est.L.3}
 |L(t,s)-L(t,s')| & \leq \kappa |s-s'|^\delta |t-s|^{-\frac{1}{2}-\delta},
 \quad \textrm{pour } |s-s'|\leq\frac{1}{2}|t-s|.
 \end{align}
 Il existe une constante $C(\delta)$ indépendante de $L$ telle que, pour tout $U \in L^2(0,1)$:
 \begin{equation}
 \label{est.L.final}
 \left| \langle LU, U \rangle \right| \leq \kappa C(\delta) |U|_{H^{-1/4}(0,1)}^2.
 \end{equation}
\end{lemme}

Les critères fournis par le Lemme~\ref{lemma.wsio} sont quantitatifs et
permettent de conclure. Ils sont liés au fait que de tels noyaux sont
faiblement singuliers. Leur dégénérescence le long de la diagonale est moins
forte que la singularité habituelle des noyaux de type Calder\'on-Zygmund.
Ainsi, les opérateurs intégraux associés font gagner de la dérivabilité.

Ce travail permet d'exhiber une quantité coercive au temps final. Au second 
ordre, le système dérive donc dans une certaine direction lors de l'application 
du contrôle $u$. Pour conclure sur le système non linéaire initial, il faut 
parvenir à estimer les restes (d'ordre cubique et plus) qui constituent la différence 
entre l'approximation de second ordre $\eta y + \eta^2 z$ et $\psi$. La difficulté est 
d'établir des estimations de ces restes ne nécessitant que la norme $H^{-5/4}$ 
du contrôle (qui est, \emph{a priori}, plus faible que les normes habituelles pour 
lesquelles le système est bien posé). Cela s'avère néanmoins possible. 
On démontre ainsi que le système, partant d'une donnée initiale nulle, dérive 
dans la direction $\rho$. Par conséquent, si on choisit une donnée initiale 
dont la projection sur $\rho$ est déjà positive, le système va continuer à 
s'éloigner de zéro ; ce qui permet de nier la contrôlabilité locale en temps 
petit.

\subsection{Commentaires}

\begin{itemize}
 \item Ce résultat est surtout intéressant pour sa méthode : lorsque le linéarisé
 ne suffit pas à conclure, on écrit le deuxième ordre, on cherche des observables
 sous la forme d'opérateurs intégraux quadratiques, on considère les opérateurs
 limites correspondant à l'asymptotique temps petit, puis on utilise la coercivité
 obtenue pour conclure sur le système initial. 
 On peut notamment envisager d'appliquer une méthode analogue 
 à des équations comme Korteweg-de-Vries ou Schrödinger dans les cas qui 
 sont encore ouverts à ce jour.
 \item Dans ce cas particulier, on démontre que l'opérateur quadratique considéré
 est coercif pour tous les contrôles. Or, il suffirait de démontrer une inégalité
 de coercivité sur l'espace des contrôles permettant de ramener le système linéarisé
  à zéro. Pour d'autres systèmes, il est probable qu'il soit important
 d'effectuer cette distinction.
 \item La question de la contrôlabilité locale de~\eqref{OQB.2} en temps long
 est ouverte. On ne sait pour le moment pas dire si
 le noyau~\eqref{OQB.7} a un signe en temps long.
 \item En dimension finie, il est possible d'expliciter la direction de dérive
 en écrivant que l'état se comporte (aux premiers ordres) 
 comme une quantité positive fois une direction fixée (dite direction de dérive).
 Ici, il n'est pas clair que cette direction existe et si elle existe, quelle
 est sa régularité (peut-être une distribution).
\end{itemize}

\section*{Conclusion}

Une question ouverte importante au sein de cette thématique est de déterminer si le Théorème~\ref{thm1} peut s'étendre en dimension infinie (l'ordre quadratique n'apporte jamais rien de bon) 
ou si au contraire il existe un système de dimension infinie pour lequel l'approximation quadratique permet d'obtenir de la contrôlabilité locale en temps petit.

\bibliographystyle{plain}
\bibliography{biblio_cmls_marbach_2017}

\end{document}